\documentclass[11pt]{article}

\usepackage{amsfonts}
\usepackage{amsmath}
\usepackage{amssymb}

\let\ds=\displaystyle
\def\R{\mathbb{R}_q^+}
\def\C{\mathbb{C}}
\def\N{\mathbb{N}}
\def\Z{\mathbb{Z}}
\def\r{\mathbb{R}}

\def\F{\mathcal{F}_{q,\v}}
\def\I{\infty}
\def\ll{\mathcal{L}_{q,1,\v}}

\def\lL{\mathcal{L}_{q,2,\v}}
\def\lp{\mathcal{L}_{q,p,\v}}

\def\v{\nu}
\newcommand{\qhyp}[5]{\mbox{}_{#1}\phi_{#2}\left(\left.\begin{array}{c}
{#3}\\{#4}\end{array}\right| q ; {#5}\right)}

\newcommand{\Hp}{\ell^2(\N)}

\newtheorem{t1}{Theorem}
\newtheorem{c1}{Corollary}
\newtheorem{d1}{Definition}

\newtheorem{l1}{Lemma}
\newtheorem{p1}{Proposition}
\newtheorem{r1}{Remark}
\newenvironment{proo}[1][Proof]{\noindent\textbf{#1.} }{\ \rule{0.5em}{0.5em}}

\begin{document}

\title{\bf Jacobi operator, $q$-difference equation and orthogonal polynomials}

\author{
Lazhar Dhaouadi
\thanks{
IPEB, 7021 Zarzouna, Bizerte,Tunisia.
E-mail : Lazhar.Dhaouadi@ipeib.rnu.tn}
\qquad \& \quad
Mohamed Jalel Atia
\thanks{
Facult\'e de Science de Gab\`es, Universit\'e de Gab\`es, Tunisie.
E-mail : jalel.atia@gmail.com}
}
\date{ }
\maketitle

\begin{abstract}
In this paper, a link between $q$-difference equations, Jacobi operators and orthogonal polynomials is given. Replacing the variable $x$ by $ q^{-n}$ in a Sturm-Liouville $q$-difference equation  we discovered the Jacobi operator. With appropriate initial conditions, the eigenfunctions of such operators are either $q$-orthogonal polynomials or the modified $q$-Bessel function and a newborn the $q$-Macdonald ones. The new Polynomial sequence we found is related to the $q$-Lommel polynomials  introduced by Koelink and other \cite{Ko,KoS,KoA}. Adapting E. C. Titchmarsh's theory, we showed the existence of a solution
square-integrable only in the complex case. As application in the real case we
gave the behavior at infinity for $q$-Macdonald's function.
Finally, we pointed out that the method described in our paper can
be generalized to study the orthogonal polynomial sequence introduced by Al-Salam and Ismail \cite{AL}.
\vspace{5mm} \\
\noindent {\it Keywords : Jacobi operator; $q$-difference equation; orthogonal polynomials; Modified $q$-Bessel functions}
\vspace{3mm}\\
\noindent {\it 2000 AMS Mathematics Subject Classification---Primary
33D15,47A05. }
\end{abstract}

\section{Introduction}
In this paper we study the relationship between three principal topics, namely

1. $q$-difference equations

2. Jacobi operators

3. Orthogonal polynomials

Each of these subjects has been intensively studied by various authors and some links have been established to joint them. For example, the theory of $q$-orthogonal polynomials was built while studying $q$-difference equations. The spectral study of Jacobi operators was started by many authors but here we showed their relationship with the $q$-difference equations. So what is new in this work?

There are two ways of seeing a $q$-difference equation:\\
- The first is to study a $q$-equation with one variable $x$
of the Sturm-Liouville kind and find power series solutions.\\
- The second rises from a simple observation, which turns out successful,
that is to replace in the $q$-difference equation the variable $x$ by $q^{-n}$ which
brings up the Jacobi operator.
With appropriate initial conditions, the eigenfunctions of such operators  are either $q$-orthogonal polynomials or special functions. This, immediately, able us to highlight the link between $q$-special functions and $q$-orthogonal polynomials. One of the consequences of this connection is  either their generating function or their measure of orthogonality.
In our case, the special functions are the modified $q$-Bessel and
among them we have a newborn the $q$-Macdonald function which we investigated its properties more closely with the great help of the $q$-Bessel Fourier analysis tools \cite{DFK}.\\
It turns out that this $q$-Macdonald function  is different from  the $q$-Macdonald function of third kind introduced by Rogov \cite{R}.
On the other hand we found a new family which is related to the $q$-Lommel polynomials introduced by Koelink and other \cite{Ko,KoS,KoA}. Thus we could determine their generating function and their orthogonality measure (which is discrete).

To enrich our paper we applied the methods of the two
classical theories in our study.

In the case of the Sturm-Liouville equation we were able to adapt E. C. Titchmarsh's theory \cite{T,S} to show the existence of a solution square-integrable only in the complex case.

In the real case we found another approach giving, as application, the behavior at infinity of the $q$-Macdonald's function.

In the case of the Jacobi operator we used the tools of moments theory to obtain some additional information on the subject.

We completed this work by noting that
the method described in our paper can
be generalized to study the orthogonal polynomial sequence introduced by Al-Salam and Ismail \cite{AL}.

\section{Preliminaries and $q$-Bessel Fourier transform}

Throughout this paper,  we consider $0<q<1$  and $\v>-1$ such that $\v\notin\N$. We adopt
the standard conventional notations of  \cite{GR}. We put $\R=\{q^n,\quad
n\in\Z\}$ and for a complex number $a$
$$
(a;q)_0=1,\quad (a;q)_n=\prod_{i=0}^{n-1}(1-aq^{i}), \quad
n=1...,+\infty,
$$
and
$$
\left[\begin{array}{c}n\\m\end{array}\right]_q=\frac{(q,q)_n}{(q,q)_m(q,q)_{n-m}}.
$$
The $q$-hypergeometric series (or basic hypergeometric) series
$_r\phi_s$ is defined by
\begin{eqnarray*}
\qhyp{r}{s}{a_1, \dots, a_r}{b_1, \dots, b_s}{z} =
\sum\limits_{k=0}
^{\infty}\frac{(a_1;q)_k(a_2;q)_k\cdots(a_r;q)_k}{(q;q)_k(b_1;q)_
k\cdots
(b_s;q)_k} \left( (-1)^kq^{\frac{1}{2}k(k-1)}\right)^{1+s-r} z^k
\label{qhyper}
\end{eqnarray*}
whenever the series converges.

\bigskip

The third Jackson $q$-Bessel function $J_{\v}$ (also called Hahn-Exton $q$-Bessel functions) is defined by the power series  \cite{S}
$$
J_{\v}(x;q)
= \frac{(q^{\v +1};q)_{\infty}}{(q;q)_{\infty}}x^{\v}
\sum_{n=0}^{\infty} \left(-1\right)^{n}
\frac{q^{\frac{n(n+1)}{2}}}{(q^{\v +1};q)_{n}(q;q)_{n}} \, x^{2n},
$$
and has the normalized form
$$
j_{\v}(x;q)=\sum_{n=0}^{\infty }(-1)^{n}\frac{q^{\frac{n(n+1)}{2}}}{
(q,q)_{n}(q^{\v +1},q)_{n}}x^{2n}.
$$
The $q$-derivative of a function $f$ is defined  for $x\neq 0$ by
$$
D_qf(x)=\frac{f(x)-f(qx)}{(1-q)x},
$$
and the $q$-Bessel operator by
$$
\Delta _{q,\v}f(x)=\frac{1}{x^{2}}\Bigg[ f(q^{-1}x)-(1+q^{2\v})f(x)+q^{2\v}f(qx)\Bigg].
$$
The $q$-Wronskian was introduced in \cite{DBF} as follows
$$
W_x(f,h)=q^{-2}(1-q)^2\Big[\Lambda_q^{-1}D_qf(x)h(x)-\Lambda_q^{-1}D_qh(x)f(x)\Big],
$$
where
$$
\Lambda_qf(x)=f(qx).
$$
The Jackson's $q$-integrals from  $0$ to $a$ and from $0$ to
$\infty$ are defined by \cite{J}
$$
\int_0^af(x)d_qx=(1-q)a\sum_{n=0}^\infty q^nf(aq^n),\quad
\int_0^\infty f(x)d_qx=(1-q)\sum_{n=-\infty}^\infty q^nf(q^n).
$$
For $1\leq p<\I$, the space $\lp$  denotes the set of real functions on
$\R$ for which
$$
\|f\|_{q,p,\v}=\left[\int_0^\I |f(x)|^px^{2\v+1}d_qx\right]^{1/p}\quad{\rm is~finite}.
$$
The $q$-Bessel Fourier transform $\mathcal{F}_{q,\v}$ is defined by
\cite{DBF,KS}
$$
\mathcal{F}_{q,\v}f(x)=c_{q,\v}\int_{0}^{\infty
}f(t)j_{\v}(xt,q^{2})t^{2\v+1}d_{q}t,\quad\forall x\in \R,
$$
where
$$
c_{q,\v}=\frac{1}{1-q}\frac{(q^{2\v+2};q^2)_\infty}{(q^{2};q^2)_\infty}.
$$
In \cite{DFK} the following results was proved:

\begin{p1}\label{pf}
For all $\v\in\r/\Z_{-}$ we have
$$
|j_{\v}(q^n,q^2)|\leq\frac{(-q^2;q^2)_\I(-q^{2\v+2};q^2)_\I}{(q^{2\v+2};q^2)_\I}\left\{
\begin{array}{c}
  1\quad\quad\quad\quad\quad\text{if}\quad n\geq 0 \\
  q^{n^2-(2\v+1)n}\quad\text{if}\quad n<0
\end{array}
\right..
$$
\end{p1}

\begin{p1}\label{pe}
Let $n,m\in\Z$ and $n\neq m$, then we have
$$
c_{q,\v}^2 \int_0^\I j_\v(q^nx,q^2)j_\v(q^mx,q^2)x^{2\v+1}d_qx=\frac{q^{-2n(\v+1)}}{1-q}\delta_{nm}.
$$
In particular
$$
\int_0^\I j_\v(q^nx,q^2)x^{2\v+1}d_qx=0.
$$
\end{p1}

\begin{t1}\label{thb}
For $f\in\ll$  we have
$$
\mathcal{F}_{q,\v}^2(f)(x)=f(x),\quad\forall x\in\R.
$$
For  $f,\F(f)\in\ll$  we have
$$
\|\mathcal{F}_{q,\v}(f)\|_{q,2,\v}=\|f\|_{q,2,\v}.
$$
\end{t1}

We end this section by the following result (see \cite{DBF}):

\begin{r1}\label{ra}
Let $f,g\in\lL$ such that $\Delta _{q,\v}f,\Delta _{q,\v}g\in\lL$. If
we have
$$
D_qf(x)=O(x^{-\v})\quad\text{and}\quad D_qg(x)=O(x^{-\v})
$$
as $x\downarrow 0$ then
$$
\langle\Delta_{q,\v}f,g\rangle=\langle f,\Delta_{q,\v}g\rangle,
$$
where $\langle .,.\rangle$ denotes the inner product on the Hilbert
space $\lL$.
\end{r1}

\section{Modified $q$-Bessel functions}

\begin{d1}We introduce  the  modified $q$-Bessel functions as follows
$$I_\v(x,q^2)=j_\v(ix,q^2),\quad i^2=-1,$$

$$\gamma_\v(x,q^2)=x^{-2\v}j_{-\v}(q^{-\v}x,q^2),$$

$$\pi_\v(x,q^2)=x^{-2\v}j_{-\v}(q^{-\v}ix,q^2),$$

and

$$
K_\v(x,q^2)=c_{q,\v}\int_0^\I\left[1+t^2\right]^{-1}j_\v(tx,q^2)t^{2\v+1}d_qt.
$$
This function can be viewed as a $q$-version of the Macdonald
function (see \cite[p.434]{W}).
\end{d1}

\begin{t1}\label{tha}
The $q$-Macdonald function $K_\v\in\ll$ and we have
\begin{equation}\label{eqa}
\mathcal{F}_{q,\v}(K_\v)(x)=\left[1+x^2\right]^{-1},\quad\forall
x\in\mathbb{R}_q^+.
\end{equation}
\end{t1}

\begin{proo}To show  $K_\v\in\ll$ we consider three cases:

\bigskip

{\bf i.} $\v>0$; {\bf ii.} $-1<\v< 0$ and {\bf iii.} $\v=0$.

\bigskip

{\bf Case i.} For $x<1$ we have
\begin{eqnarray*}
x^{2\v+2}|K_\v(x,q^2)|&\leq&c_{q,\v}\int_0^\I
\left[1+\frac{t^2}{x^2}\right]^{-1}|j_\v(t,q^2)|t^{2\v+1}d_qt\\
&\leq&\left\{c_{q,\v}\int_0^\I
\left[1+\frac{x^2}{t^2}\right]^{-1}|j_\v(t,q^2)|t^{2\v-1}d_qt\right\}x^2\\
&\leq&\left\{c_{q,\v}\int_0^\I |j_\v(t,q^2)|t^{2\v-1}d_qt\right\}x^2
\end{eqnarray*}
and then $\ds\int_0^1|K_\v(x,q^2)|x^{2\v+1}d_qx<\I.$ On the other hand, for
$x>1$ we have
\begin{eqnarray*}
x^2K_\v(x,q^2)&=&c_{q,\v}\int_0^\I
\left[1+t^2\right]^{-1}x^2j_\v(xt,q^2)t^{2\v+1}d_qt\\
&=&-c_{q,\v}\int_0^\I
\left[1+t^2\right]^{-1}\Delta_{q,\v}j_\v(xt,q^2)t^{2\v+1}d_qt\\
&=&-c_{q,\v}\int_0^\I
\Delta_{q,\v}\left[1+t^2\right]^{-1}j_\v(xt,q^2)t^{2\v+1}d_qt\quad \text{(from Remark \ref{ra})}\\
&=&-c_{q,\v}\int_0^\I
u(t)j_\v(xt,q^2)t^{2\v+1}d_qt\\
\end{eqnarray*}
where $u(t)=\Delta_{q,\v}\left[1+t^2\right]^{-1}$ is a bounded function
on $\R$. Multiplying $x^2K_\nu(x,q^2)$ by $x^{2\nu+2}$
and making an upper bound for the quantity we obtain
\begin{eqnarray*}
x^2\Big[x^{2\v+2}|K_\v(x,q^2)|\Big]&\leq&c_{q,\v}\int_0^\I
|u(t/x)||j_\v(t,q^2)|t^{2\v+1}d_qt\\
&\leq&\left\{c_{q,\v}\|u\|_{q,\I}\int_0^\I
|j_\v(t,q^2)|t^{2\v+1}d_qt\right\},
\end{eqnarray*}
then $\ds\int_1^\I|K_\v(x,q^2)|x^{2\v+1}d_qx<\I$.

\bigskip

{\bf Case ii.}  In fact $\ds\lim_{x\to 0}K_\v(x,q^2)$ exists and then
$\ds\int_0^1|K_\v(x)|x^{2\v+1}d_qx<\I.$ Using the same method as in case
i, we prove that $\ds\int_1^\I|K_\v(x,q^2)|x^{2\v+1}d_qx<\I.$

\bigskip

{\bf Case iii.} For $x<1$ we have
\begin{eqnarray*}
x^{2}|K_0(x,q^2)|&\leq&c_{q,0}\int_0^\I
\left[1+\frac{t^2}{x^2}\right]^{-1}|j_0(t,q^2)|td_qt\\
&\leq&\left\{c_{q,0}\int_0^\I
\left(\frac{t}{x}\right)\left[1+\frac{t^2}{x^2}\right]^{-1}|j_0(t,q^2)|d_qt\right\}x\\
&\leq&\left\{c_{q,0}\int_0^\I |j_0(t,q^2)|d_qt\right\}x
\end{eqnarray*}
and then $\ds\int_0^1|K_0(x,q^2)|x^{2v+1}d_qx<\I.$ The proof is identical
to case ii.

\bigskip

This proves that $K_\v\in\ll$  and then (\ref{eqa}) is a simple  consequence of the inversion formula in
Theorem \ref{thb}
\end{proo}

\begin{p1}\label{pc}
Let $f$ and $h$ two linearly independent solutions of the following $q$-difference equation
$$
\Delta_{q,\v}y(x)=\pm~~ y(x).
$$
Then there exists a constant $c(f,h)\neq 0$ such that
$$
x^{2\v}\Big[f(x)h(qx)-f(qx)h(x)\Big]=c(f,g),\quad\forall x\in\R.
$$
\end{p1}

\begin{proo}
With
$$
D_q\Big[y\mapsto y^{2\v+1}W_y(f,h)\Big](x)=\Big[\Delta_{q,\v}f(x)h(x)-h(x)\Delta_{q,\v}f(x)\Big]x^{2\v+1}=0,
$$
there exists a constant $c$ such that
$$
q^{-2}(1-q)^2x^{2\v+1}\Big[\Lambda_q^{-1}D_qf(x)h(x)-\Lambda_q^{-1}D_qh(x)f(x)\Big]=c,\quad\forall x\in\R.
$$
To completes the proof, it suffice to use the definition of the operator $D_q$.
\end{proo}

\begin{p1}\label{pa}
The functions
$x\mapsto j_\v(\lambda x,q^2)$ and $x\mapsto \gamma_\v(\lambda x,q^2)$  are two  linearly independent solutions of the following equation
$$
\Delta_{q,\v}y(x)=-\lambda^2 y(x).
$$
In addition
$$
c(j_\v,\gamma_\v)=(q^{-2\v}-1).
$$
\end{p1}

\begin{proo}For all $\alpha\in\r$ such that $\alpha\neq -1,-2,\ldots$ we have \cite{DFK}
$$
\Delta _{q,\alpha }j_{\alpha }(\lambda x,q^{2})=-\lambda ^{2}j_{\alpha}(\lambda x,q^{2}).
$$
Let $\mu=q^{-\v}\lambda$, then we have
\begin{eqnarray*}
\Delta _{q,\v}\gamma _{\v}(\lambda x,q^{2}) &=&x^{-2v}\frac{
q^{2v}j_{-\v}(q^{-1-\v}\lambda x,q^{2})-(1+q^{2\v})j_{-\v}(q^{-\v}\lambda
x,q^{2})+j_{-\v}(q^{1-\v}\lambda x,q^{2})}{x^{2}} \\
&=&q^{2\v}x^{-2\v}\frac{j_{-\v}(q^{-1-\v}\lambda
x,q^{2})-(1+q^{-2\v})j_{-\v}(q^{-\v}\lambda
x,q^{2})+q^{-2\v}j_{-\v}(q^{1-\v}\lambda x,q^{2})}{x^{2}} \\
&=&q^{2\v}x^{-2\v}\frac{j_{-\v}(q^{-1}\mu x,q^{2})-(1+q^{-2\v})j_{-\v}(\mu
x,q^{2})+q^{-2\v}j_{-\v}(q\mu x,q^{2})}{x^{2}} \\
&=&q^{2\v}x^{-2\v}\left[ -\mu ^{2}j_{-\v}(\mu x,q^{2})\right] =-\lambda
^{2}\gamma _{\v}(\lambda x,q^{2}).
\end{eqnarray*}
On the other hand we have
$$
c(j_\v,\gamma _\v)=q^{-2v}j_\v(x,q^2)j_{-\v}(q^{1-\v}x,q^2)-j_\v(qx,q^2)j_{-\v}(q^{-\v}x,q^2).
$$
When $x\to 0$ we obtain the result.
\end{proo}

\begin{p1}\label{pb}
The functions
$x\mapsto I_\v(\lambda x,q^2)$ and $x\mapsto K_\v(\lambda x,q^2)$  are two  linearly independent solutions of the following equation
\begin{equation}\label{eqh}
\Delta_{q,\v}y(x)=\lambda^2 y(x).
\end{equation}
\end{p1}

\begin{proo}Let $g_{\lambda}(x)=\lambda^{2(\v+1)}K_\v(\lambda x,q^2)$ then
$$
g_{\lambda}(x)=c_{q,\v}\int_0^\I\left[1+\frac{t^2}{\lambda^2}\right]^{-1}j_\v(tx,q^2)t^{2\v+1}d_qt.
$$
\begin{eqnarray*}
\left[1-\frac{\Delta_{q,\v}}{\lambda^2}\right]g_{\lambda}(x)&=&c_{q,\v}\int_0^\I \left[1+\frac{t^2}{\lambda^2}\right]
^{-1}\left[1-\frac{\Delta_{q,\v}}{\lambda^2}\right]j_\v(tx,q^2)t^{2\v+1}d_qt\\
&=&c_{q,\v}\int_0^\I j_\v(tx,q^2)t^{2\v+1}d_qt=0.
\end{eqnarray*}
Note that
$$
\left[1-\frac{\Delta_{q,\v}}{\lambda^2}\right]j_\v(tx,q^2)=\left[1+\frac{t^2}{\lambda^2}\right]j_\v(tx,q^2).
$$
The function $K_\v\in\ll$ but $I_\v\notin\ll$. Hence, we conclude that they provide two linearly independent solutions.
\end{proo}

\begin{l1}Let $\lambda\in\C$ such that $\lambda\notin\R\cup q^\v\R$ then we have
\begin{equation}\label{eqg}
\lim_{n\rightarrow \infty }q^{2\v n}\frac{j_{-\v}(q^{-n-\v}\lambda ,q^{2})}{
j_{\v}(q^{-n}\lambda ,q^{2})}=\frac{(q^{2\v+2},q^{2})_{\infty }}{(q^{-2\v+2},q^{2})_{\infty }}\frac{(q^{2\v}\lambda ^{-2},q^{2})_{\infty
}(q^{2-2\v}\lambda ^{2},q^{2})_{\infty }}{(\lambda ^{-2},q^{2})_{\infty
}(q^{2}\lambda ^{2},q^{2})_{\infty }}.
\end{equation}
\end{l1}

\begin{proo}
Let $x\in\C^*\backslash\R$ then we have the following asymptotic expansion as $|x|\to\I$
$$
j_{\alpha}(x,q^{2})\sim \frac{(x^{2}q^{2},q^{2})_{\infty }}{(q^{2\alpha
+2},q^{2})_{\infty }},\quad\forall\alpha\in\r\backslash\Z_{-}.
$$
If $n\to\I$ we obtain
$$
j_\v(q^{-n}\lambda ,q^{2})\sim \frac{(q^{2-2n}\lambda
^{2},q^{2})_{\infty }}{(q^{2\v +2},q^{2})_{\infty }},\quad
\forall\lambda \notin\R.
$$
On the other hand we have
$$
(q^{2-2n}\lambda ^{2},q^{2})_{\infty }=(-1)^{n}\lambda
^{2n}q^{-n(n-1)}(\lambda ^{-2},q^{2})_{n}\left( q^{2}\lambda
^{2},q^{2}\right) _{\infty }.
$$
Hence, when $n\to\I$
$$
j_{\v}(q^{-n}\lambda ,q^{2})\sim \frac{(-1)^{n}\lambda
^{2n}q^{-n(n-1)}(\lambda ^{-2},q^{2})_{n}(q^{2}\lambda ^{2},q^{2})_{\infty }
}{(q^{2\v+2},q^{2})_{\infty }},\quad\forall\lambda \notin\R,
$$
and
$$
j_{-\v}(q^{-n-\v}\lambda ,q^{2})\sim \frac{(-1)^{n}q^{-2\v n}\lambda
^{2n}q^{-n(n-1)}(q^{2\v}\lambda ^{-2},q^{2})_{n}(q^{2-2\v}\lambda
^{2},q^{2})_{\infty }}{(q^{-2\v+2},q^{2})_{\infty }},\quad
\forall\lambda\notin q^\v\R.
$$
This implies
$$
q^{2\v n}\frac{j_{-v}(q^{-n-\v}\lambda ,q^{2})}{j_{\v}(q^{-n}\lambda ,q^{2})}
\sim \frac{(q^{2\v+2},q^{2})_{\infty }}{(q^{-2\v+2},q^{2})_{\infty }}\frac{(q^{2\v}\lambda ^{-2},q^{2})_{n}(q^{2-2\v}\lambda
^{2},q^{2})_{\infty }}{(\lambda ^{-2},q^{2})_{n}(q^{2}\lambda
^{2},q^{2})_{\infty }},\quad\forall\lambda\notin\R\cup q^\v\R.
$$
Hence,
$$
\lim_{n\rightarrow \infty }q^{2vn}\frac{j_{-\v}(q^{-n-\v}\lambda ,q^{2})}{
j_{\v}(q^{-n}\lambda ,q^{2})}=\frac{(q^{2\v+2},q^{2})_{\infty }}{(q^{-2\v+2},q^{2})_{\infty }}\frac{(q^{2\v}\lambda ^{-2},q^{2})_{\infty
}(q^{2-2\v}\lambda ^{2},q^{2})_{\infty }}{(\lambda ^{-2},q^{2})_{\infty
}(q^{2}\lambda ^{2},q^{2})_{\infty }},\quad\forall\lambda\notin\R\cup q^v\R.
$$
This completes the proof.
\end{proo}

\begin{p1}We have
\begin{equation}\label{eqi}
K_\v(x,q^2)=\alpha_\v\Big[\pi_\v(x,q^2)-\beta_\v I_\v(x,q^2)\Big]
\end{equation}
where
$$
\beta_\v=\frac{(q^{2v+2},q^{2})_{\infty }}{(q^{-2\v+2},q^{2})_{\infty }}\frac{(-q^{2\v},q^{2})_{\infty
}(-q^{2-2\v},q^{2})_{\infty }}{(-1,q^{2})_{\infty
}(-q^{2},q^{2})_{\infty }},
$$
and
$$
\alpha_\v=\left\{\begin{array}{c}
  \frac{(q^2,q^2)_\I}{(q^{2\v},q^2)_\I}\quad\quad\quad\quad\quad\quad\quad\quad\quad{\rm if}\quad \v\geq 0 \\
  \\
  -\frac{c_{q,\v}}{\beta_\v}\ds\int_0^\I [1+t^2]^{-1}t^{2\v+1}d_qt\quad{\rm if}\quad \v<0
\end{array}\right..
$$
\end{p1}

\begin{proo}
The functions
$x\mapsto I_\v(\lambda x,q^2)$ and $x\mapsto \pi_\v(\lambda x,q^2)$  are two  linearly independent solutions of (\ref{eqh}).
Then there exist two constants $\alpha_\v$ and $\beta_\v$ such that (\ref{eqi}) hold true. Now we can write
$$
K_\v(q^{-n},q^2)=\alpha_\v\Big[\frac{\pi_\v(q^{-n},q^2)}{I_\v(q^{-n},q^2)}-\beta_\v\Big]I_\v(q^{-n},q^2).
$$
On the other hand
$$
\lim_{n\to\I}I_\v(q^{-n},q^2)=+\I
$$
Using Theorem \ref{tha} we see that
$$
\lim_{n\to\I}K_\v(q^{-n},q^2)=0.
$$
Then  it is necessary that
$$
\lim_{n\to\I}\Big[\frac{\pi_\v(q^{-n},q^2)}{I_\v(q^{-n},q^2)}-\beta_\v\Big]=0.
$$
Formula (\ref{eqg}) with $\lambda=i$ lead to the result. To evaluate $\alpha_\v$ we consider two cases:

If $\v>0$ then
$$
\alpha_\v=\lim_{x\to 0}x^{2\v}K_\v(x,q^2)=c_{q,\v}\int_0^\I j_\v(t,q^2)t^{2\v-1}d_qt.
$$
Using an identity established in \cite{KS}
\begin{multline*}
\int_{0}^{\infty} t^{-\lambda } J_{\theta}(q^{n}t;q^{2}) J_{\mu }(q^{m}t;q^{2}) \, d_{q}t  \\
= (1-q)q^{n(\lambda -1)+(m-n)\mu }
\frac{(q^{1+\lambda+\theta-\mu},q^{2\mu +2};q^{2})_{\infty}}
{(q^{1-\lambda +\theta +\mu},q^{2};q^{2})_{\infty}}
\\
\times {_{2}\phi_{1}}
\left(\,
\begin{matrix}
q^{1-\lambda +\mu +\theta },q^{1-\lambda+\mu-\theta} \\
q^{2\mu +2}
\end{matrix}
\,\middle|\,
q^{2};q^{2m-2n+1+\lambda+\theta-\mu} \,\right),
\end{multline*}
valid for $n,m\in\Z$ and $\theta,\mu,\lambda\in\r$ such that
$$1-\lambda+\theta+\mu>0,\quad 2m-2n+1+\lambda+\theta-\mu>0,$$
we conclude that $(\theta=\v,\mu=0,\lambda=1-\v,n=0,m\to\I)$
$$
\alpha_\v=\frac{(q^2,q^2)_\I}{(q^{2\v},q^2)_\I}.
$$
If $-1<\v<0$ then
$$
\alpha_\v=-\frac{1}{\beta_\v}\lim_{x\to 0}K_\v(x,q^2)=-\frac{c_{q,\v}}{\beta_\v}\int_0^\I \Big[1+t^2\Big]^{-1}t^{2\v+1}d_qt.
$$
\end{proo}

\begin{c1} As a direct consequence, we have
$$
c(I_\v,K_\v)=\alpha_\v(q^{-2\v}-1).
$$
\end{c1}

\section{Jacobi operator and $q$-difference equation}

A tridiagonal matrix of the form
$$
J =\begin{pmatrix}  b_0 & a_0 & 0   & 0   & 0 & 0 & \ldots\cr
             a_0 & b_1 & a_1 & 0   & 0 & 0 & \ldots\cr
             0  & a_1 & b_2 & a_2  & 0 & 0 & \ldots\cr
             0  & 0   & a_2 & b_3& a_3 & 0 & \ldots\cr
             \vdots&  &\ddots &\ddots &\ddots &\ddots&
   \end{pmatrix}
$$
is a Jacobi operator, or an infinite Jacobi matrix,
if $b_i\in\r$ and $a_i>0$.

We consider $J$ as an operator defined on the
Hilbert space $\Hp$. So with respect
to the standard orthonormal basis $\{ e_k\}_{k\in\N}$ of $ \Hp$
the Jacobi operator is defined as
$$
J\, e_k = \begin{cases} a_k\, e_{k+1} + b_k\, e_k +
a_{k-1} \, e_{k-1}, & k\geq 1, \\
a_0\, e_1 + b_0\, e_0, & k=0. \end{cases}
$$
Associated with the Jacobi matrix $J$ there is a
three-term recurrence relation
\begin{equation}\label{eqb}
a_nP_{n+1}(\lambda)+ b_nP_n(\lambda)+ a_{n-1}P_{n-1}(\lambda)=\lambda P_n(\lambda),\quad n\geq0:
\end{equation}
If we take the initial conditions $P_{-1}(x)=0$ and $P_0(x)=1$ we obtain a sequence of orthogonal
polynomials $(P_n)$. We assume that $P_n$ is of degree $n$ with positive leading coefficient.

\bigskip

The polynomials of the second kind $(Q_n)$ are generated
by the three-term recurrence relation (\ref{eqb}) but with initial conditions $Q_{-1}(x)=-1$ and $Q_0(x)=0$.
Consequently, $(P_n)$ and $(Q_n)$ are linearly independent solutions of (\ref{eqb}) and together they span the solution
space. Notice that $Q_n(x)$ is a polynomial of degree $n-1$.
The orthonormal polynomials $(P_n)$ and the polynomials of the second kind $(Q_n)$ play a crucial role for the theory
of the moment problem.

\bigskip

If $(P_n)$ satisfy the above three-term recurrence relation (including the initial conditions) for some real sequences
$(a_n)$ and $(b_n)$ with $a_n>0$, then it follows by Favard's theorem that there exists a positive measure $\mu$ on
$\r$ such that the polynomials $(P_n)$ are orthonormal with respect to $\mu$.

The recurrence coefficients in (\ref{eqb}) contain useful information about the moment problem. Carleman
proved in 1922 that the moment problem is determinate if
\begin{equation}\label{eqf}
\sum_{n=0}^\I\frac{1}{a_n}=\I.
\end{equation}
If $\mu$ is determinate, then $\sum_{n=0}^\I |P_n(z)|^2=\I$ for all $z\in\C\backslash\r$ and all $z\in\r$ except for at most countably many points $z$ where $\mu(\{z\})>0$, in the case of which
$$
\mu(\{z\})=\left(\sum_{n=0}^\I P_n^2(z)\right)^{-1}.
$$
If $\mu$ is indeterminate, then
\begin{equation}\label{e5}
\sum_{n=0}^\I |P_n(z)|^2<\I,\quad \sum_{n=0}^\I |Q_n(z)|^2<\I,\quad\forall z\in\C.
\end{equation}

There exists close relationship between $q$-difference and Jacobi operators. Let us describe briefly this relationship.

\bigskip

Given an eigenfunction $f$ of the following  $q$-difference operator
$$
H_{q,\v}f(x)=\Delta_{q,\v}f(x)-V(x)f(x)=-\lambda f(x),\quad\forall x\in\R,
$$
then we have
$$
\frac{1}{x^2}f(q^{-1}x)-\Bigg[\frac{1+q^{2\v}}{x^2}+V(x)\Bigg]f(x)+\frac{q^{2\v}}{x^2}f(qx)=-\lambda f(x).
$$
Assume that $f(x)=x^{-(\v+1)}g(x)$, then
\begin{equation}\label{ee}
\frac{q^{\v+1}}{x^2}g(q^{-1}x)-\Bigg[\frac{1+q^{2\v}}{x^2}+V(x)\Bigg]g(x)+\frac{q^{\v-1}}{x^2}g(qx)=-\lambda g(x).
\end{equation}
If we replace $x$ by $q^{-n}$ and $g(q^{-n})$ by $g_n$ in (\ref{ee}) then we obtain
$$
Jg_n=a_ng_{n+1}+b_ng_n+a_{n-1}g_{n-1}=-\lambda g_n,
$$
where
\begin{equation}\label{e6}
a_n=q^{2n+\v+1},\quad b_n=-q^{2n}(1+q^{2\v})-V(q^{-n}).
\end{equation}
Note that:

. $\{g_n\}$ is an eigenfunction of the Jacobi operator $J$ if and only if $f$ is an eigenfunction of the  $q$-difference operator $H_{q,\v}$.

. $\{g_n\}\in\ell^2(\N)$ if and only if $\ds\int_a^\I |f(x)|^2x^{2\v+1}d_qx<\I$ for some constant $a\in\R$.

. $\{g_n\}\in\ell^2(\Z)$ if and only if $\ds\int_0^\I |f(x)|^2x^{2\v+1}d_qx<\I$.

\begin{t1}\label{t1}
We assume that  $\lambda$ is a complex number other than real values. There exists a unique solution $f$ (up to constant factor) of the $q$-difference equation
\begin{equation}\label{e2}
\Delta_{q,\v}f(x)=[V(x)-\lambda]f(x),
\end{equation}
satisfying
$$
\int_a^\I |f(x)|^2x^{2\v+1}d_qx<\I,
$$
for some constant $a\in\R$.
\end{t1}

\begin{proo}
If $F(x)=F(x,\lambda)$  and $G(x)=G(x,\lambda')$ satisfy (\ref{e2}) then,
\begin{align*}
& (\lambda^{\prime}-\lambda)\int_{a}^{b}F(x)G(x)x^{2\v+1}d_qx\\
& =\int_{a}^{b}\Big[F(x)\Big\{  V(x)G(x)-\Delta_{q,\v}G(x)\Big\} -G(x)\Big\{
V(x)F(x)-\Delta_{q,\v}F(x)\Big\}\Big]x^{2\v+1}d_qx\\
& =-\int_{a}^{b}\Big[F(x)\Delta_{q,\v}G(x)-G(x)\Delta_{q,\v}F(x)\Big]
x^{2\v+1}d_qx=a^{2\v+1}W_{a}(F,G)-b^{2\v+1}W_{b}(F,G).
\end{align*}
where $W_{a}(F,G)$ is the $q$-wronskian.

If $\lambda=\mu+i\eta$, $\lambda^{\prime}=\overline{\lambda}$ and
$G=\overline{F}$, then we get
$$
2\nu\int_{a}^{b}\left\vert F(x)\right\vert ^{2}x^{2\v+1}d_qx=ia^{2\v+1}W_{a}(F,\overline
{F})-ib^{2\v+1}W_{b}(F,\overline{F}).
$$
In the sequel we assume that $\eta>0$.

\bigskip

Now let $\phi(x)=\phi(x,\lambda)$ and $\theta(x)=\theta(x,\lambda)$ be two
solutions of (\ref{e2}) such that
$$
\left\{
\begin{array}[c]{c}
\phi(a)=\epsilon\sin\alpha,\quad\Lambda_{q}^{-1}D_{q}\phi(a)=-\epsilon\cos\alpha\\
\theta(a)=\epsilon\cos\alpha,\quad\Lambda_{q}^{-1}D_{q}\theta(a)=\epsilon\sin\alpha
\end{array}
\right.  ,
$$
where $\alpha$ is real and
$$
\epsilon=\frac{1}{\sqrt{q^{-1}(1-q)^2a^{2\v+1}}}.
$$
It follows that
$$
b^{2\v+1}W_{b}(\phi,\theta)=a^{2\v+1}W_{a}(\phi,\theta)=\sin^{2}\alpha+\cos^{2}\alpha=1,
$$
and
$$
a^{2\v+1}W_{a}(\phi,\overline{\phi})=a^{2\v+1}W_{a}(\theta,\overline{\theta})=0.
$$
Therefore
$$
a^{2\v+1}W_{a}(\theta+l\phi,\overline{\theta}+\overline{l}\overline{\phi}
)=l-\overline{l}=2i\operatorname{Im}l.
$$
Let $z$ be a complex variable such that
$$
l=l(\lambda)=-\frac{\theta(b)z+\Lambda_{q}^{-1}D_{q}\theta(b)}{\phi(b)z+\Lambda_{q}^{-1}D_{q}\phi(b)}.
$$
If $z$ describes the real axis $D=\{z\in\C,\quad\operatorname{Im}z=0\}$ then $l$ describes a circle $C_b$ in the complex plane.

\bigskip

Here $l=\I$ correspond to $z=-\ds\frac{\Lambda_{q}^{-1}D_{q}\phi(b)}{\phi(b)}$.
The center $l_b$ of $C_{b}$ correspond to the conjugate, therefore
$$
l_b=-\frac{W_{b}(\theta,\overline{\phi})}{W_{b}(\phi,\overline{\phi})}.
$$
Also
\begin{align*}
&\operatorname{Im}\left[-\frac{\Lambda_{q}^{-1}D_{q}\phi(b)}{\phi(b)}\right]
=-\frac{1}{2}i\left[\frac{\Lambda_{q}^{-1}D_{q}\overline{\phi}(b)}{\overline{\phi}(b)}-\frac{\phi_{q}^{-1}D_{q}\phi(b)}{\phi(b)}\right]\\
&=-\frac{1}{2}i\frac{W_{b}(\phi,\overline{\phi})}{|\phi(b)|^2}
=\frac{\eta}{b^{2\v+1}|\phi(b)|^2}\int_a^b|\phi(x)|^2x^{2\v+1}d_qx>0.\\
\end{align*}
The exterior of $C_b$ correspond to $\operatorname{Im}z>0$. Now $l$ is inside $C_b$  if $\operatorname{Im}z<0$

i.e. if
$$
i\left[  -\frac{l\Lambda_{q}^{-1}D_{q}\phi(b)+\Lambda_{q}^{-1}D_{q}\theta
(b)}{l\phi(b)+\theta(b)}+\frac{\overline{l}\Lambda_{q}^{-1}D_{q}\overline
{\phi}(b)+\Lambda_{q}^{-1}D_{q}\overline{\theta}(b)}{l\overline{\phi
}(b)+\overline{\theta}(b)}\right]  >0,
$$
equivalently
$$
i\left[  \left\vert l\right\vert ^{2}W_{b}(\phi,\overline{\phi})+lW_{b}%
(\phi,\overline{\theta})+\overline{l}W_{b}(\theta,\overline{\phi}%
)+W_{b}(\theta,\overline{\theta})\right]  >0,
$$
that is
$$
iW_{b}(\theta+l\phi,\overline{\theta}+\overline{l}\overline{\phi})>0,
$$
equivalently
$$
\int_{a}^{b}\left\vert \theta+l\phi\right\vert^{2}x^{2\v+1}d_qx<-\frac{\operatorname{Im}l}{\eta}.
$$
Since $-\ds\frac{\Lambda_{q}^{-1}D_{q}\theta(b)}{\Lambda_{q}^{-1}D_{q}\phi(b)}$ is on
$C_{b}$ (for $z=0$) the radius $r_{b}$ of $C_{b}$ is
$$
r_{b}=\left\vert \frac{\Lambda_{q}^{-1}D_{q}\theta(b)}{\Lambda_{q}^{-1}
D_{q}\phi(b)}-\frac{W_{b}(\theta,\overline{\phi})}{W_{b}(\phi,\overline{\phi})}\right\vert =\left\vert \frac{b^{2\v+1}W_{b}(\theta,\phi)}{b^{2\v+1}W_{b}(\phi,\overline
{\phi})}\right\vert =\frac{1}{2\eta\int_{a}^{b}\left\vert \phi(x)\right\vert
^{2}x^{2\v+1}d_{q}x}.
$$
If $b<b^{\prime}$ then $C_{b^{\prime}}$ includes $C_{b}$ ($r_{b'}<r_b$). As $b\rightarrow\infty$, the
circle $C_{b}$ converge either to a limit-circle or to a limit-point.

\bigskip

If $m=m(\lambda)$ is the limit-point, or any point on the limit-circle,
$$
\int_{a}^{b}\left\vert \theta+m\phi\right\vert ^{2}x^{2\v+1}d_{q}x<-\frac
{\operatorname{Im}m}{\eta},
$$
for all values of $b$. Hence
$$
\int_{a}^{\infty}\left\vert \theta+m\phi\right\vert ^{2}x^{2\v+1}d_{q}x<-\frac
{\operatorname{Im}m}{\eta}.
$$
Note that if $\operatorname{Im}(\lambda)=\eta<0$ then it suffices to consider the following $q$-difference equation
$$
\Delta_{q,v}\overline{f}(x)=\Big[\overline{V}(x)-\overline{\lambda}\Big]\overline{f}(x).
$$
Using (\ref{e6}) and (\ref{eqf}) we see that the moment problem is determinate.  By (\ref{e5}) we conclude that when the moment problem is indeterminate we are in the limit-circle at infinity. So in our case we are in the limit-point at infinity.
\end{proo}

\bigskip

In the following we assume that  $\lambda$ is a real value and  there exist $a\in\R$ and $\delta>0$ such that
\begin{equation}\label{e4}
V(x)-\lambda >\delta,\quad \forall x\in\R,\quad x\geq a.
\end{equation}
We denote by
$$
f_{n}=f(q^{-n}),\quad \v_{n}=q^{-\v}\frac{f_{n}}{f_{n-1}},\quad {\rm if}\quad f_{n-1}\neq 0.
$$

\begin{l1}\label{l1}
Assume that $f(x)\neq 0$ for $x$ big enough and
\begin{equation}\label{e3}
\int_a^\I |f(x)|^2x^{2\v+1}d_qx<\I,
\end{equation}
then when $n$ is big enough we have
$$
|\v_{n}| \leq 1.
$$
\end{l1}

\begin{proo}
\bigskip In fact if $a=q^{-n_0}$ then
$$
\sum_{n=n_0}^\I q^{-2(\v+1)n}| f(q^{-n})|^2 <\I,
$$
which implies
$$
\limsup_{n\rightarrow \I }\Bigg[\frac{q^{-(\v+1)(n+1)}| f(q^{-n-1})|
}{q^{-(\v+1)n}| f(q^{-n})| }\Bigg]\leq 1.
$$
When $n$ is big enough we have
$$
\Bigg[\frac{q^{-(\v+1)(n+1)}| f(q^{-n-1})|
}{q^{-(\v+1)n}| f(q^{-n})| }\Bigg]\leq \frac{1}{q},
$$
and the result follows.
\end{proo}

\begin{t1}\label{t2}
Let $f$ be a solution of the $q$-difference equation
$$
\Delta _{q,\v}f(x)=\left[ V(x)-\lambda \right] f(x),\quad\forall
x\in\R.
$$
Then either $\ds\lim_{x\to\I}f(x)=\pm \I$ or $f$ satisfies the conditions of Lemma \ref{l1} and then  when
$n$ big enough we have:

\begin{description}
  \item[i.] $\v_{n}>0$ and $\ds\lim_{n\rightarrow \I }\v_{n}=0.$
  \item[ii.] There exist $c,\sigma>0$ such that
$$
|f_{n}| <\sigma c^{n}q^{n^{2}}.
$$
\end{description}
\end{t1}

\begin{proo}We have
\begin{eqnarray*}
q(1-q)^{-2}x^{2\v+1}\left[ V(x)-\lambda \right] f(x)
&=&q(1-q)^{-2}x^{2\v+1}\Delta _{q,\v}f(x) \\
&=&D_{q}\Big[y\mapsto y^{2v+1}D_{q}f(q^{-1}y)\Big](x).
\end{eqnarray*}
Let $a\in\R$ such that
$$
V(x)-\lambda >\delta,\quad \forall x\geq a.
$$
Assume that $f(a)>0.$

\begin{description}
\item[$\star$] If $D_{q}f(a)>0$ then
$$
D_{q}f(q^{-n}a)>D_{q}f(a)>0,
$$
which implies
$$
f(q^{-n}a)>(q^{-n}a)D_{q}f(a)>0.
$$
This gives
$$
\lim_{n\to +\I}f(q^{-n})=+\I .
$$
\item[$\star$] If $D_{q}f(a)<0$\ then there are three possibilities:
\begin{description}
\item[]
\item[$\bullet$] The function $f(x)$ changes sign, it can do so at most once and then
goes to $-\I .$ (replace $f$ by $-f$ and use previous argument).

\item[$\bullet$] The function $f(x)$ remains positive and  $D_{q}f(x)$ change
sign then $f(x)$ goes to $+\I .$

\item[$\bullet$] The function $f(x)$ remains positive and  $D_{q}f(x)$  remains negative for $x\geq a$, then
$f$ decreases steadily, and then goes to a positive limit. Hence
$D_{q}f(x)\rightarrow 0$. Moreover
\begin{eqnarray*}
\int_{a}^{b}\Delta _{q,\v}f(x)x^{2\v+1}d_{q}x &=&\int_{a}^{b}\left[
V(x)-\lambda \right] f(x)x^{2\v+1}d_{q}x \\
&=&q^{-1}(1-q)^{2}\Big[b^{2\v+1}D_{q}f(q^{-1}b)-a^{2\v+1}D_{q}f(q^{-1}a)\Big]
\\
&<&-q^{-1}(1-q)^{2}a^{2\v+1}D_{q}f(q^{-1}a),\quad \forall b>a.
\end{eqnarray*}
Hence
$$
\delta\int_{a}^{\I}f(x)x^{2\v+1}d_{q}x<\int_{a}^{\I}\left[V(x)-\lambda \right] f(x)x^{2\v+1}d_{q}x<-q^{-1}(1-q)^{2}a^{2\v+1}D_{q}f(q^{-1}a),
$$
which proves that $f$ satisfies (\ref{e3}).
\end{description}
\end{description}
Similar arguments hold if $f(a)<0$.

\bigskip

Assume that $f$ satisfies the conditions of Lemma \ref{l1}. We have
$$
f(q^{-1}x)-(1+q^{2\v})f(x)+q^{2\v}f(qx)=x^{2}\left[ q(x)-\lambda \right] f(x).
$$
Let
$$
q_{n}=V(q^{-n}),
$$
then we have
$$
f_{n+1}-(1+q^{2\v})f_{n}+q^{2\v}f_{n-1}=q^{-2n}\left[ q_{n}-\lambda \right]
f_{n}
$$
$$
q^{-\v}\frac{f_{n+1}}{f_{n}}+q^{\v}\frac{f_{n-1}}{f_{n}}=q^{-\v}\Big[ q^{-2n}
\left[q_{n}-\lambda \right] +(1+q^{2\v})\Big] .
$$
Set
$$
A_{n}=q^{-\v}\Bigg[ q^{-2n}\left[ q_{n}-\lambda \right] +\left(1+q^{2\v}\right)\Bigg] ,
$$
we get
$$
\v_{n+1}=A_{n}-\frac{1}{\v_{n}}.
$$
From Lemma \ref{l1} we see that
$$
\left|A_{n}-\frac{1}{\v_{n}}\right| =|\v_{n+1}| \leq 1.
$$
As $V(x)$ satisfies the condition (\ref{e4}), we obtain
$$
\lim_{n\rightarrow \I }A_{n}=+\I .
$$
Hence
$$
\left\{
\begin{array}{c}
\v_{n}>0 \\
\ds\lim_{n\rightarrow \I }\v_{n}=0
\end{array}
\right..
$$
On the other hand
$$
\left|q^{-2n-\v}\left[ q_{n}-\lambda \right] -\frac{1}{\v_{n}}\right|
<\left( q^\v+q^{-\v}+1\right) .
$$
Hence, as $n$ tends to $+\I$
$$
\left|1-\frac{1}{\v_{n}}\frac{q^{2n+\v}}{\left[ q_{n}-\lambda \right]}
\right| <\frac{\left( q^\v+q^{-\v}+1\right) q^{\v}}{\delta
}q^{2n}\rightarrow 0.
$$
Then when $n$ is big enough we have
\begin{align*}
& \frac{1}{\v_{n}}\frac{q^{2n+v}}{\left[ q_{n}-\lambda
\right] }>
\frac{1}{2} \\
& \Rightarrow 2\frac{q^{\v}}{\delta }q^{2n}>\v_{n} \\
& \Rightarrow |f_{n}| <\left[ 2\frac{q^{2\v}}{\delta } \right]
q^{2n}|f_{n-1}|,
\end{align*}
which leads to the result.
\end{proo}

\begin{c1}
The $q$-Mackdonal  function satisfies the following properties when $n$ is big enough
\begin{description}
  \item[i.] The sequence $K_{\v}(q^{-n},q^2)$ have a constant sign and
$$
\lim_{n\rightarrow \I }\frac{K_{\v}(q^{-n},q^2)}{K_{\v}(q^{-n+1},q^2)}=0.
$$
  \item[ii.] There exist $c,\sigma>0$ such that
$$
|K_{\v}(q^{-n},q^2)| <\sigma c^{n}q^{n^{2}}.
$$
\end{description}
\end{c1}

\begin{proo}
From Theorem \ref{tha} we see that the $q$-Mackdonal function satisfies the condition of Lemma \ref{l1}. Then we apply Theorem \ref{t2} with $V=0$ and $\lambda=-1$. This give the result.
\end{proo}

\section{Orthogonal polynomial}

Let $V=0$ in (\ref{e6}) and  consider the three-term recurrence relation
\begin{equation}\label{eqc}
 q^{2n+\v+1}P_{n+1}(\lambda)-q^{2n}(1+q^{2\v})P_n(\lambda)+ q^{2n+\v-1}P_{n-1}(\lambda)=\lambda P_n(\lambda).
\end{equation}
with the initial condition  $P_{-1}=0$ and $P_{0}=1$.
Using (\ref{eqf}) we see that the moment problem is determinate and so,
there exists a unique positive measure $\mu$ on $\r$
such that the polynomials $(P_n)$ are orthonormal with respect to $\mu$.

\begin{r1}
The $q$-Lommel polynomials $r_n(.,w,q)$  are given as follows \cite{Ko}
$$
\left\vert w\right\vert ^{-1}q^{-\frac{1}{2}
-n}r_{n+1}(x,w,q)+q^{-n}(1+w^{-2})r_{n}(x,w,q)+\left\vert w\right\vert
^{-1}q^{\frac{1}{2}-n}r_{n-1}(x,q,w)=xr_{n}(x,w,q)
$$
with initial conditions
$$
r_{-1}(x,w,q)=0,\quad r_{0}(x,w,q)=1
$$
where
\begin{equation}\label{e7}
0<q<1,\quad |w|>1.
\end{equation}
So it is easy to show that if $w=q^{-\v}$ then we obtain
$$
P_{n}(x)=r_{n}(x,w,q^{-2}).
$$
Here $P_n(x)$ are the polynomials satisfying (\ref{eqc}) and are not a particular case of the $q$-Lommel
polynomials because they do not satisfy the conditions (\ref{e7}). We can consider $P_n(x)$ as the $q^{\it -1}$-Lommel
polynomials.\bigskip

Moreover if we replace $x$ by $q^n$ and $g(q^n)$ by $g_n$ in (\ref{ee}) then the recurrence relation (\ref{eqc}) becomes
\begin{equation*}
 q^{-2n+\v-1}P_{n+1}(\lambda)-q^{-2n}(1+q^{2\v})P_n(\lambda)+ q^{2n+\v+1}P_{n-1}(\lambda)=\lambda P_n(\lambda).
\end{equation*}
So with the choice $w=q^{-v}$ and the condition $v>0$ we obtain $P_{n}(x)=r_{n}(x,w,q^2)$
which are a particular case of the $q$-Lommel polynomials. The polynomial sequence \eqref{eqc} is not $q$-classical, in fact, from \cite[p.65-67]{KM}, we can see that the recurrence coefficients are neither classical nor shifted-classical.
\end{r1}

\begin{p1}\label{pd}
Let $f$ and $h$ be two linearly independent solutions of the following $q$-difference equation
\begin{equation}\label{eqe}
\Delta_{q,\v}y(x)=\lambda y(x).
\end{equation}
Then
$$
P_{n}(\lambda)=q^{-n(\v+1)}\left[\frac{f(q^{-n})h(q)-h(q^{-n})f(q)}{f(1)h(q)-h(1)f(q)}\right].
$$
\end{p1}

\begin{proo} The polynomials $(P_n)$ and $(Q_n)$ are linearly independent solutions of (\ref{eqc}) and together they form the space
of solutions. Then there exist two functions $\theta$ and $m$ such that
$$
g_{n}=\theta (\lambda)P_{n}(\lambda)+m(\lambda)Q_{n}(\lambda).
$$
We have
$$
f(q^{-n})=q^{n(\v+1)}g_{n}.
$$
This implies
$$
f(q^{-n})=q^{n(\v+1)}\Big[\theta (\lambda)P_{n}(\lambda)+m(\lambda)Q_{n}(\lambda)\Big].
$$
For $n=-1$ we get
$$
m(\lambda)=-q^{(\v+1)}f(q),
$$
and for $n=0$
$$
\theta(\lambda)=f(1).
$$
Then
\begin{equation}\label{eqd}
f(q^{-n})=q^{n(\v+1)}\left[f(1)P_{n}(\lambda)-q^{(\v+1)}f(q)Q_{n}(\lambda)\right].
\end{equation}
Also for the function $h$ we have
$$
h(q^{-n})=q^{n(\v+1)}\left[h(1)P_{n}(\lambda)-q^{(\v+1)}h(q)Q_{n}(\lambda)\right],
$$
which give the result.
\end{proo}

\begin{c1}\label{ca}
We have the following identity
\begin{description}
  \item[a)] $j_\v(q^{-n}\lambda,q^2)=q^{n(\v+1)}\Big[j_\v(\lambda,q^2)P_{n}(-\lambda^2)-q^{(\v+1)}j_\v(q\lambda,q^2)Q_{n}(-\lambda^2)\Big].$
  \item[b)] $\gamma_\v(q^{-n}\lambda,q^2)=q^{n(\v+1)}\Big[\gamma_\v(\lambda,q^2)P_{n}(-\lambda^2)-q^{(\v+1)}\gamma_\v(q\lambda,q^2)Q_{n}(-\lambda^2)\Big].$
  \item[c)] $I_\v(q^{-n}\lambda,q^2)=q^{n(\v+1)}\Big[I_\v(\lambda,q^2)P_{n}(\lambda^2)-q^{(\v+1)}I_\v(q\lambda)Q_{n}(\lambda^2)\Big].$
  \item[d)] $K_\v(q^{-n}\lambda,q^2)=q^{n(\v+1)}\Big[K_\v(\lambda,q^2)P_{n}(\lambda^2)-q^{(\v+1)}K_\v(q\lambda,q^2)Q_{n}(\lambda^2)\Big].$
  \item[e)] $\ds P_{n}(-\lambda^2)=q^{-n(\v+1)}\frac{\lambda^{2\v}}{(q^{-2\v}-1)}\Big[j_\v(q^{-n}\lambda,q^2)\gamma_\v(q\lambda,q^2)-\gamma_\v(q^{-n}\lambda,q^2)j_\v(q\lambda,q^2)\Big].$
  \item[f)] $\ds P_{n}(\lambda^2)=q^{-n(\v+1)}\frac{\lambda^{2\v}}{\alpha_\v(q^{-2\v}-1)}\Big[I_\v(q^{-n}\lambda,q^2)K_\v(q\lambda,q^2)-K_\v(q^{-n}\lambda,q^2)I_\v(q\lambda,q^2)\Big].$
\end{description}
\end{c1}

\begin{proo}The identity a) and b) are  simple consequences of (\ref{eqd}) if we replace $\lambda$ by $-\lambda^2$ in equation (\ref{eqe}). For the identity c) and d) we replace $\lambda$ by $\lambda^2$. To prove e) and f) we use Propositions \ref{pc} and \ref{pd}.
\end{proo}

\begin{c1}\label{cb}
We have
$$
P_n(\lambda)=\frac{q^{-n(\v+1)}}{q^{-2\v}-1}\sum_{m=0}^{n}a_{m,n}\lambda^{m},
$$
where
\begin{align*}
a_{m,n}&=\frac{q^{m(m+1)-2m\v}}{(q^{2},q^{2})_{m}(q^{-2\v+2},q^{2})_{m}}\times\\
&\left\{
q^{2(m-\v)}\text{ }_{2}\phi _{1}\left[\left.
\begin{array}{c}
q^{-2m},q^{-2m+2\v} \\
q^{2v+2}
\end{array}
\right\vert q^{2},q^{2(\v-n)}\right] -q^{2n(m+\v)}\text{ }_{2}\phi _{1}\left[
\left.
\begin{array}{c}
q^{-2m},q^{-2m+2\v} \\
q^{2\v+2}
\end{array}
\right\vert q^{2},q^{4+2(\v-n)}\right] \right\}.
\end{align*}
\end{c1}

\begin{proo} Using the identities e) in Corollary \ref{ca} we obtain
\begin{align*}
P_{n}(-\lambda^{2})  & =\frac{q^{-n(\v+1)}}{c(j_{\v},\gamma_{\v})}\lambda
^{2v}\Big[j_{\v}(q^{-n}\lambda,q^2)\gamma_{\v}(q\lambda,q^2)-j_{\v}(q\lambda,q^2
)\gamma_{\v}(q^{-n}\lambda,q^2)\Big]  \\
& =\frac{q^{-n(\v+1)}}{c(j_{\v},\gamma_{v})}\Big[q^{-2\v}j_{\v}(q^{-n}
\lambda,q^2)j_{-\v}(q^{1-\v}\lambda,q^2)-q^{2n\v}j_{\v}(q\lambda,q^2)j_{-\v}(q^{-n-\v}
\lambda,q^2)\Big].
\end{align*}
Taking into account formula \cite[(6.4.4)]{S} for the
product of two Hahn-Exton $q$-Bessel functions
$$
j_{\alpha }(ax,q^{2})j_{\beta }(bx,q^{2})=\sum_{m=0}^{\infty }(-1)^{m}\frac{%
q^{m(m+1)}b^{2m}}{(q^{2},q^{2})_{m}(q^{2\beta +2},q^{2})}\text{ }_{2}\phi
_{1}\left[ \left.
\begin{array}{c}
q^{-2m},q^{-2m-2\beta } \\
q^{2\alpha +2}
\end{array}
\right\vert q^{2},\left( q^{\alpha +\beta +1}\frac{a}{b}\right) ^{2}\right]
x^{2m}
$$
we obtain
$$
q^{-2\v}j_{\v}(q^{-n}x,q^{2})j_{-\v}(q^{1-\v}x,q^{2})=\sum_{m=0}^{\infty
}(-1)^{m}\frac{q^{m(m+1)+2m(1-\v)-2\v}}{(q^{2},q^{2})_{m}(q^{-2\v+2},q^{2})_{m}}%
\text{ }_{2}\phi _{1}\left[ \left.
\begin{array}{c}
q^{-2m},q^{-2m+2\v} \\
q^{2\v+2}
\end{array}
\right\vert q^{2},q^{2(\v-n)}\right] x^{2m}
$$
and
$$
q^{2n\v}j_{\v}(qx,q^{2})j_{-\v}(q^{n-\v}x,q^{2})=\sum_{m=0}^{\infty }(-1)^{m}
\frac{q^{m(m+1)+2m(n-\v)+2n\v}}{(q^{2},q^{2})_{m}(q^{-2\v+2},q^{2})_{m}}\text{ }
_{2}\phi _{1}\left[ \left.
\begin{array}{c}
q^{-2m},q^{-2m+2\v} \\
q^{2\v+2}
\end{array}
\right\vert q^{2},q^{4+2(\v-n)}\right] x^{2m}
$$
which give the expression of the coefficients $a_{m,n}$.

\bigskip

Note that $a_{m,n}=0$ if $m>n$.
\end{proo}

\begin{p1}\label{ph}
The generating function for the  polynomials $P_n$ is
$$
\sum_{n=0}^{\infty }t^{n}P_{n}(x)=\frac{1}{(1-q^{\v-1}t)(1-q^{-(\v+1)}t)}\text{ }_{2}\phi
_{2}\left[ \left.
\begin{array}{c}
q^{2},0 \\
q^{1+\v}t,q^{1-\v}t
\end{array}
\right\vert q^{2};-q^{-(1+\v)}xt\right],\quad q^{1-\v}t~~{\rm and}~~q^{-(1+\v)}t\notin\Z_{-}
$$
As a consequence we have
$$
P_{n}(x)=q^{-n(1+\v)}\sum_{m=0}^{n}q^{m(m-1)}\left[
\begin{array}{c}
n \\
m%
\end{array}
\right] _{q^{2}}\text{ }_{2}\phi _{1}\left[ \left.
\begin{array}{c}
q^{2m+2},q^{2(m-n)} \\
q^{-2n}
\end{array}
\right\vert q^{2};q^{2(\v-m)}\right] x^{m}.
$$
\end{p1}

\begin{proo}
As in \cite[p. 13]{KoS}  we will derive a generating function for the  polynomials $P_n$. Multiply
(\ref{eqc}) with $t^{n+1}$ and sum from $n=0$ to $\I$ we have
$$
q^{\v-1}\sum_{n=0}^{\infty }\left( q^{2}t\right)
^{n+1}P_{n+1}(x)-(1+q^{2\v})t\sum_{n=0}^{\infty }\left( q^{2}t\right)
^{n}P_{n}(x)+q^{\v+1}t^{2}\sum_{n=0}^{\infty }\left( q^{2}t\right)
^{n-1}P_{n-1}(x)=xt\sum_{n=0}^{\infty }\left( q^{2}t\right) ^{n}P_{n}(x).
$$
Let us introduce the generating function
$$
G(x,t)=\sum_{n=0}^{\infty }t^{n}P_{n}(x)
$$
then one gets
$$
G(x,q^{2}t)=\frac{1+q^{1-\v}xtG(x,t)}{\left( 1-q^{1+\v}t\right) \left(
1-q^{1-\v}t\right)}.
$$
With $G(x,0)=1$ we obtain
\begin{eqnarray*}
G\left( x,t\right) &=&\sum_{m=0}^{\infty }\frac{q^{m(m-1)-(1+\v)m}}{\left(
q^{\v-1}t,q^{2}\right) _{m+1}\left( q^{-(1+\v)}t,q^{2}\right) _{m+1}}\left(
xt\right) ^{m} \\
&=&\frac{1}{(1-q^{\v-1}t)(1-q^{-(\v+1)}t)}\text{ }_{2}\phi _{2}\left[ \left.
\begin{array}{c}
q^{2},0 \\
q^{1+\v}t,q^{1-\v}t
\end{array}
\right\vert q^{2};-q^{-(1+v)}xt\right].
\end{eqnarray*}
Using the $q$-binomial theorem one obtains
$$
\frac{1}{\left( q^{\v-1}t,q^{2}\right) _{m+1}}=\text{ }_{1}\phi _{0}\left[
\left.
\begin{array}{c}
q^{2m+2} \\
-
\end{array}
\right\vert q^{2};q^{\v-1}t\right],\quad |q^{\v-1}t|<1.
$$

$$
\frac{1}{\left( q^{-(\v+1)}t,q^{2}\right) _{m+1}}=\text{ }_{1}\phi _{0}\left[
\left.
\begin{array}{c}
q^{2m+2} \\
-
\end{array}
\right\vert q^{2};q^{-(\v+1)}t\right],\quad |q^{-(\v+1)}t|<1.
$$
This yields
\begin{eqnarray*}
G\left( x,t\right) &=&\sum_{m,j,k=0}^{\infty }q^{m(m-1)-(1+\v)m+j(\v-1)-k(\v+1)}
\frac{\left( q^{2m+2},q^{2}\right) _{j}\left( q^{2m+2},q^{2}\right) _{k}}{
\left( q^{2},q^{2}\right) _{j}\left( q^{2},q^{2}\right) _{k}}x^{m}t^{m+j+k}
\\
&=&\sum_{n=0}^{\infty }q^{-n(1+\v)}\left( \sum_{m+j+k=n}q^{m(m-1)+2j\v}\frac{
\left( q^{2m+2},q^{2}\right) _{j}\left( q^{2m+2},q^{2}\right) _{k}}{\left(
q^{2},q^{2}\right) _{j}\left( q^{2},q^{2}\right) _{k}}x^{m}\right) t^{n}.
\end{eqnarray*}
Then we obtain the explicit representation
\begin{eqnarray*}
P_{n}(x) &=&q^{-n(1+\v)}\left( \sum_{m+j+k=n}q^{m(m-1)+2j\v}\frac{\left(
q^{2m+2},q^{2}\right) _{j}\left( q^{2m+2},q^{2}\right) _{k}}{\left(
q^{2},q^{2}\right) _{j}\left( q^{2},q^{2}\right) _{k}}x^{m}\right)  \\
&=&q^{-n(1+\v)}\sum_{m=0}^{n}q^{m(m-1)}\left( \sum_{j+k=n-m}\frac{\left(
q^{2m+2},q^{2}\right) _{j}\left( q^{2m+2},q^{2}\right) _{k}}{\left(
q^{2},q^{2}\right) _{j}\left( q^{2},q^{2}\right) _{k}}q^{2j\v}\right) x^{m} \\
&=&q^{-n(1+\v)}\sum_{m=0}^{n}q^{m(m-1)}\left( \sum_{j=0}^{n-m}\frac{\left(
q^{2m+2},q^{2}\right) _{j}\left( q^{2m+2},q^{2}\right) _{n-m-j}}{\left(
q^{2},q^{2}\right) _{j}\left( q^{2},q^{2}\right) _{n-m-j}}q^{2j\v}\right)
x^{m}.
\end{eqnarray*}
the following identities
$$
\left( q^{2},q^{2}\right) _{n-m-j}=(-1)^{j}\frac{\left( q^{2},q^{2}\right)
_{n-m}}{\left( q^{2(m-n)},q^{2}\right) _{j}}q^{j(j-1)-2(n-m)j}
$$

$$
(q^{2m+2},q^{2})_{n-m-j}=(-1)^{j}\frac{\left( q^{2m+2},q^{2}\right) _{n-m}}{
\left( q^{-2n},q^{2}\right) _{j}}q^{j(j-1)-2(n-m)j-2mj}
$$
leads to
$$
\frac{\left( q^{2m+2},q^{2}\right) _{n-m-j}}{\left( q^{2},q^{2}\right)
_{n-m-j}}=\frac{\left( q^{2m+2},q^{2}\right) _{n-m}}{\left(
q^{2},q^{2}\right) _{n-m}}\frac{\left( q^{2(m-n)},q^{2}\right) _{j}}{\left(
q^{-2n},q^{2}\right) _{j}}q^{-2mj},
$$
which gives
\begin{eqnarray*}
P_{n}(x) &=&q^{-n(1+\v)}\sum_{m=0}^{n}q^{m(m-1)}\frac{\left( q^{2m+2},q^{2}\right)
_{n-m}}{\left( q^{2},q^{2}\right) _{n-m}}\left( \sum_{j=0}^{n-m}\frac{\left(
q^{2m+2},q^{2}\right) _{j}\left( q^{2(m-n)},q^{2}\right) _{j}}{\left(
q^{2},q^{2}\right) _{j}\left( q^{-2n},q^{2}\right) _{j}}q^{2(\v-m)j}\right)
x^{m} \\
&=&q^{-n(1+\v)}\sum_{m=0}^{n}q^{m(m-1)}\left[
\begin{array}{c}
n \\
m
\end{array}
\right] _{q^{2}}\left( \sum_{j=0}^{n-m}\frac{\left( q^{2m+2},q^{2}\right)
_{j}\left( q^{2(m-n)},q^{2}\right) _{j}}{\left( q^{2},q^{2}\right)
_{j}\left( q^{-2n},q^{2}\right) _{j}}q^{2(\v-m)j}\right) x^{m} \\
&=&q^{-n(1+\v)}\sum_{m=0}^{n}q^{m(m-1)}\left[
\begin{array}{c}
n \\
m
\end{array}
\right] _{q^{2}}\text{ }_{2}\phi _{1}\left[ \left.
\begin{array}{c}
q^{2m+2},q^{2(m-n)} \\
q^{-2n}
\end{array}
\right\vert q^{2};q^{2(\v-m)}\right] x^{m}.
\end{eqnarray*}
\end{proo}

\begin{t1}\label{thc}
The Stieltjes transform of the orthogonality measure $\mu$ for the orthogonality polynomials $P_n$ is given for all $z\notin{\rm Supp}(\mu)$ by
\begin{multline*}
\int_{\r}\frac{d\mu(t)}{z-t}=\frac{1}{q^{(\v+1)}j_{\v}(q\sqrt{z} ,q^{2})}\times\\
\left[ j_{\v}(\sqrt{z} ,q^{2})-\frac{
(q^{-2\v}-1)(z^{-1},q^{2})_{\infty }(q^{2}z,q^{2})_{\infty}}{q^{2\v}j_{-\v}(q^{1-\v}\sqrt{z} ,q^{2})(z^{-1},q^{2})_{\infty }(q^{2}z,q^{2})_{\infty }-\sigma_\v j_{\v}(q\sqrt{z} ,q^{2})(q^{2\v}z^{-1},q^{2})_{\infty}(q^{2-2\v}z,q^{2})_{\infty}}\right],
\end{multline*}
where
$$
\sigma_\v=\frac{(q^{2\v+2},q^{2})_{\infty }}{(q^{-2\v+2},q^{2})_{\infty }}.
$$
\end{t1}

\begin{proo}
The Stieltjes transform can be obtained from \cite[thm. 2.4]{AI}
$$
\int_{\r}\frac{d\mu(t)}{z-t}=\lim_{n\to\I}\frac{Q_n(z)}{P_n(z)},\quad\forall z\notin{\rm Supp}(\mu).
$$
Using the identity a) in Corollary \ref{ca}
$$
\frac{j_{\v}(q^{-n}\lambda ,q^{2})}{q^{n(\v+1)}P_{n}(-\lambda ^{2})}
=j_{\v}(\lambda ,q^{2})\text{{}}-q^{(\v+1)}j_{\v}(q\lambda ,q^{2})\frac{
Q_{n}(-\lambda ^{2})}{P_{n}(-\lambda ^{2})},
$$
and the identity e)
$$
\frac{q^{n(\v+1)}P_{n}(-\lambda ^{2})}{j_{\v}(q^{-n}\lambda ,q^{2})}=\frac{1}{
q^{-2\v}-1}\left[ q^{2v}j_{-\v}(q^{1-\v}\lambda ,q^{2})-j_{\v}(q\lambda
,q^{2})\left( q^{2n\v}\frac{j_{-\v}(q^{-n-\v}\lambda ,q^{2})}{
j_{\v}(q^{-n}\lambda ,q^{2})}\right) \right]
$$
we obtain
\begin{eqnarray*}
\frac{Q_{n}(-\lambda ^{2})}{P_{n}(-\lambda ^{2})} &=&\frac{1}{
q^{(\v+1)}j_{\v}(q\lambda ,q^{2})}\left[ j_{\v}(\lambda ,q^{2})-\frac{
j_{\v}(q^{-n}\lambda ,q^{2})}{q^{n(v+1)}P_{n}(-\lambda ^{2})}\right] \\
&=&\frac{1}{q^{(\v+1)}j_{\v}(q\lambda ,q^{2})}\left[ j_{\v}(\lambda ,q^{2})-
\frac{q^{-2\v}-1}{q^{2\v}j_{-\v}(q^{1-\v}\lambda ,q^{2})-j_{\v}(q\lambda
,q^{2})\left( q^{2n\v}\frac{j_{-\v}(q^{-n-\v}\lambda ,q^{2})}{
j_{\v}(q^{-n}\lambda ,q^{2})}\right) }\right].
\end{eqnarray*}
Hence, when we use (\ref{eqg}) we obtain
$$
\lim_{n\rightarrow \infty }\frac{Q_{n}(-\lambda ^{2})}{P_{n}(-\lambda ^{2})}=
\frac{1}{q^{(\v+1)}j_{\v}(q\lambda ,q^{2})}\left[ j_{\v}(\lambda ,q^{2})-\frac{
q^{-2\v}-1}{q^{2\v}j_{-\v}(q^{1-\v}\lambda ,q^{2})-\sigma_\v j_{\v}(q\lambda ,q^{2})\frac{
(q^{2\v}\lambda ^{-2},q^{2})_{\infty }(q^{2-2v}\lambda ^{2},q^{2})_{\infty }}{
(\lambda ^{-2},q^{2})_{\infty }(q^{2}\lambda ^{2},q^{2})_{\infty }}}\right].
$$
Now, if $\lambda\in\mathbb{R}_{q}$ and using Proposition \ref{pf} we have
$$
\lim_{n\rightarrow \infty
}\left\vert q^{2n\v}\frac{j_{-\v}(q^{-n-v}\lambda ,q^{2})}{j_{\v}(q^{-n}\lambda
,q^{2})}\right\vert =+\infty,
$$
which implies
$$
\lim_{n\rightarrow \infty }\frac{Q_{n}(-\lambda ^{2})}{
P_{n}(-\lambda ^{2})}=\frac{j_{\v}(\lambda ,q^{2})}{q^{(\v+1)}j_{\v}(q\lambda
,q^{2})}.
$$
If $\lambda\in q^\v\mathbb{R}_{q}$ then
$$
\lim_{n\rightarrow \infty
}\left\vert q^{2n\v}\frac{j_{-\v}(q^{-n-\v}\lambda ,q^{2})}{j_{\v}(q^{-n}\lambda,q^{2})}\right\vert =0.
$$
Hence we obtain
$$
\lim_{n\rightarrow \infty }\frac{Q_{n}(-\lambda ^{2})}{
P_{n}(-\lambda ^{2})}=\frac{1}{q^{(\v+1)}j_{\v}(q\lambda ,q^{2})}\left[
j_{\v}(\lambda ,q^{2})-\frac{q^{-2\v}-1}{q^{2\v}j_{-\v}(q^{1-\v}\lambda ,q^{2})}
\right].
$$
To complete the proof (case $z=\lambda^2$) we replace $\lambda$ by $i\lambda$.
\end{proo}

\begin{c1}\label{cc}
The orthogonality measure $\mu$ is supported on a denumerable discrete set
$$
\{t_k=q^{-2}(j_k^\v)^2,\quad k\in\N\}
$$
where $j_k^\v$ are the positive simple zero of the function $j_\v(x,q^2)$ numbered increasing. The positive mass
$A_k$ of $\mu$ at $t_k$ is given by
$$
A_k=\frac{2j_k^\v}{q^{(v+3)}j'_{\v}(j_k^\v,q^{2})}\times\left[ j_{\v}(q^{-1}j_k^\v,q^{2})-\frac{
(q^{-2\v}-1)}{q^{2\v}j_{-v}(q^{-\v}j_k^\v ,q^{2})}\right].
$$
\end{c1}

\begin{proo}
From Theorem \ref{thc} we deuce that the orthogonality measure $\mu$ is supported on a denumerable discrete set of the zeros $\{t_k\}$ of $j_{\v}(q\sqrt{z} ,q^{2})$. So let $\mu$ have mass $A_k$ at the point $\{t_k\}$ then we have for all $z\neq t_k$ \cite[(5.15)]{KoA}
\begin{multline*}
\sum_{k=0}^\I\frac{A_k}{z-t_k}=\frac{1}{q^{(\v+1)}j_{\v}(q\sqrt{z} ,q^{2})}\times\\
\left[ j_{\v}(\sqrt{z} ,q^{2})-\frac{
(q^{-2\v}-1)(z^{-1},q^{2})_{\infty }(q^{2}z,q^{2})_{\infty}}{q^{2\v}j_{-\v}(q^{1-\v}\sqrt{z} ,q^{2})(z^{-1},q^{2})_{\infty }(q^{2}z,q^{2})_{\infty }-\sigma_\v j_{\v}(q\sqrt{z} ,q^{2})(q^{2\v}z^{-1},q^{2})_{\infty}(q^{2-2\v}z,q^{2})_{\infty}}\right].
\end{multline*}
The zeros of $j_{v}(q\sqrt{z} ,q^{2})$ correspond precisely to the poles $t_k$ of the left hand side and $A_k$ is the residue at this point.
\end{proo}

\begin{r1}We can generalize our method to study other families of orthogonal polynomials. For instance
let $V(t)=-\frac{1+q^{2\v}}{t^{2}}$ in (\ref{e6}), then the corresponding orthogonal polynomials $P_n(x)$ satisfy
$$
a_{n}P_{n+1}(x)+b_{n}P_{n}(x)+a_{n-1}P_{n-1}(x)=xP_{n}(x),
$$
where
$$
a_{n}=q^{2n+1+\v},\quad b_{n}=0.
$$
The family of orthogonal polynomial
$$
R_n(x)=q^{n-n^2}P_n\left(q^{-(1+\v)}x\right)
$$
satisfying the following recurrence relation
$$
R_{n+1}(x)+q^{4(n-1)}R_{n-1}(x)=xR_{n}(x)
$$
with  initial conditions
$$
R_{-1}(x)=0,\quad R_{0}(x)=1.
$$
The polynomials $R_{n}(x)$ have the remarkable properties that $R_{n}(q^{2n}x)$
are also orthogonal and are the only modulo a constant (see the reference  \cite{AL} for more detail). Obviously, using the
same idea exposed above we can find the corresponding orthogonality measure and generating function and probably
some other properties.
\end{r1}

\end{document}